\RequirePackage{ifpdf}
\ifpdf % We are running pdfTeX in pdf mode
\documentclass[pdftex]{sigma}
\else
\documentclass{sigma}
\fi

\numberwithin{equation}{section} \numberwithin{theorem}{section}
\numberwithin{proposition}{section}
\numberwithin{corollary}{section} \numberwithin{remark}{section}
\numberwithin{definition}{section} \numberwithin{example}{section}

\begin{document}

\renewcommand{\thefootnote}{$\star$}

\renewcommand{\PaperNumber}{067}

\FirstPageHeading

\ShortArticleName{Inversion Formulas for the Dunkl Intertwining
Operator}

\ArticleName{Inversion Formulas for the Dunkl Intertwining\\
Operator and Its Dual on Spaces of Functions\\ and
Distributions\footnote{This paper is a contribution to the Special
Issue on Dunkl Operators and Related Topics. The full collection
is available at
\href{http://www.emis.de/journals/SIGMA/Dunkl_operators.html}{http://www.emis.de/journals/SIGMA/Dunkl\_{}operators.html}}}

\Author{Khalifa TRIM\`ECHE}

\AuthorNameForHeading{K. Trim\`eche}

\Address{Faculty of Sciences of Tunis, Department of Mathematics,
1060 Tunis, Tunisia}
\Email{\href{mailto:Khlifa.trimeche@fst.rnu.tn}{Khlifa.trimeche@fst.rnu.tn}}

\ArticleDates{Received May 13, 2008, in f\/inal form September 16,
2008; Published online September 29, 2008}

\Abstract{In this paper  we prove inversion formulas for the Dunkl
intertwining opera\-tor~$V_k$ and for its dual ${}^tV_k$
 and we
deduce the expression of the representing distributions of the
inverse operators $V^{-1}_k$ and ${}^tV_k^{-1}$, and we give some
applications.}


\Classification{33C80; 43A32; 44A35; 51F15}

\section{Introduction}\label{sec1}

We consider the dif\/ferential-dif\/ference operators $T_j$, $j =
1, 2,\dots,d$,
 on $\mathbb{R}^d$ associated to a root system $R$ and
 a multiplicity function $k$, introduced by
C.F.~Dunkl in \cite{3} and called the Dunkl ope\-rators in the
literature. These operators are very important in pure
 mathematics and in~physics.
They provide a useful tool in the study of
 special functions related to
root systems \cite{4,6, 2}. Moreover the commutative algebra
generated by these operators has been used in the study of certain
exactly solvable models of quantum mechanics, namely the
Calogero--Sutherland--Moser models, which deal with systems of
identical particles in a one dimensional space (see \cite{8,11,
12}).

C.F.~Dunkl  proved in \cite{4} that there exists a unique
isomorphism $V_k$ from the space of homogeneous polynomials
$\mathcal{P}_n$ on $\mathbb{R}^d$ of degree $n$ onto itself
satisfying the transmutation relations
\begin{gather}\label{eq1.1}
T_jV_k = V_k \frac{\partial}{\partial x_j},\qquad j = 1,
2,\dots,d,
\end{gather}
 and
 \begin{gather}
 V_k(1) = 1.\label{eq1.2}
 \end{gather}
 This operator is called the Dunkl intertwining
operator. It has been extended to an isomorphism from
$\mathcal{E}(\mathbb{R}^d)$ (the space of $C^\infty$-functions on
$\mathbb{R}^d)$ onto itself satisfying the relations \eqref{eq1.1}
and \eqref{eq1.2} (see~\cite{15}).

The operator $V_k$ possesses the integral  representation
\begin{gather}
\label{eq1.3} \forall\; x \in \mathbb{R}^d,\qquad V_k(f)(x) =
\int_{\mathbb{R}^d}f(y)d\mu_x(y), \qquad f \in
\mathcal{E}(\mathbb{R}^d),
\end{gather} where $\mu_x$ is a
probability measure on $\mathbb{R}^d$ with \mbox{support} in the
closed ball $B(0, \|x\|)$ of center $0$ and radius $\|x\|$
(see~\cite{14, 15}).

We have shown in \cite{15} that for each $x \in \mathbb{R}^d$,
there exists a unique distribution~$\eta_x$ in
$\mathcal{E}'(\mathbb{R}^d)$ (the space of distributions on
$\mathbb{R}^d$ of compact \mbox{support}) with \mbox{support} in
$B(0, \|x\|)$ such that
\begin{gather}
V_k^{-1}(f)(x) = \langle \eta_x,f\rangle, \qquad f \in \mathcal{E}
(\mathbb{R}^d).\label{eq1.4}
\end{gather}

We have studied also in \cite{15} the transposed operator
${}^tV_k$ of the operator $V_k$, satisfying for~$f$ in
$\mathcal{S}(\mathbb{R}^d)$ (the space of $C^\infty$-functions on
$\mathbb{R}^d$ which are rapidly decreasing
 together with their
derivatives) and $g$ in $\mathcal{E}(\mathbb{R}^d)$, the relation
\begin{gather*}
\int_{\mathbb{R}^d}{}^tV_k(f)(y)g(y)dy =
\int_{\mathbb{R}^d}V_k(g)(x) f(x) \omega_k(x)dx, %\label{eq1.5}
\end{gather*}
where $\omega_k$ is a positive weight function on $\mathbb{R}^d$
which will be def\/ined in the following section. It has the
integral representation
\begin{gather}
\forall\; y \in \mathbb{R}^d,\qquad {}^tV_k(f)(y) =
\int_{\mathbb{R}^d}f(x)d\nu_y(x), \label{eq1.6}
\end{gather}
 where $\nu_y$ is
a positive measure on $\mathbb{R}^d$ with \mbox{support} in the
set $\{x \in \mathbb{R}^d ;\;\|x\| \geq \|y\|\}$. This operator is
called the dual Dunkl intertwining operator.

We have proved in \cite{15} that the operator ${}^tV_k$ is an
isomorphism from $\mathcal{D}(\mathbb{R}^d)$ (the space of
$C^\infty$-functions on $\mathbb{R}^d$  with compact support)
(resp.~$\mathcal{S} (\mathbb{R}^d))$ onto itself, satisfying the
transmutation relations
\begin{gather*}
\forall\; y \in \mathbb{R}^d,\qquad {}^tV_k(T_jf)(y) =
\frac{\partial}{\partial y_j} {}^tV_k(f)(y),\qquad j = 1,
2,\dots,d.%\label{eq1.7}
\end{gather*} Moreover for each $y \in \mathbb{R}^d$,
there exists a unique distribution $Z_y$ in
$\mathcal{S}'(\mathbb{R}^d)$ (the space of tempered distributions
on $\mathbb{R}^d)$ with \mbox{support} in the set $\{x \in
\mathbb{R}^d ; \|x\| \geq \|y\|\}$
 such that
 \begin{gather}
 {}^tV_k^{-1}(f)(y) = \langle
Z_y, f\rangle ,\qquad f \in \mathcal{S} (
\mathbb{R}^d).\label{eq1.8}
\end{gather} Using the operator $V_k$,
C.F.~Dunkl has def\/ined in \cite{5} the Dunkl kernel $K$ by
\begin{gather}
\forall\; x \in \mathbb{R}^d,\qquad  \forall\; z \in
\mathbb{C}^d,\qquad K(x, -iz) = V_k(e^{-i\langle \cdot ,
z\rangle})(x).\label{eq1.9}
\end{gather} Using
this kernel C.F.~Dunkl has introduced in \cite{5} a Fourier
transform $\mathcal{F}_D$ called the Dunkl transform.

In this paper we establish the following inversion formulas for
the operators $V_k$ and ${}^tV_k$:
\begin{gather}
\forall\; x \in \mathbb{R}^d,\qquad
 V^{-1}_k(f)(x) = P{}^tV_k(f)(x),\qquad
f \in \mathcal{S}( \mathbb{R}^d),\label{eq1.10}\\
\forall\; x \in \mathbb{R}^d, \qquad {}^tV_k^{-1}(f)(x) =
V_k(P(f))(x),\qquad f \in \mathcal{S}(\mathbb{R}^d),\nonumber
% \label{eq1.11}
 \end{gather}
where $P$ is a pseudo-dif\/ferential operator on $\mathbb{R}^d$.

When the multiplicity function takes integer values, the formula
\eqref{eq1.10} can also be written in the form
\begin{gather*}
\forall\; x \in \mathbb{R}^d,\qquad V^{-1}_k(f)(x)
= {}^tV_k(Q(f))(x), \qquad f \in \mathcal{S}(\mathbb{R}^d), %\label{eq1.12}
\end{gather*}
where $Q$ is a dif\/ferential-dif\/ference operator.

Also we give another expression of the operator ${}^tV_k^{-1}$ on
the space $\mathcal{E}'( \mathbb{R}^d)$.
 From these
relations we deduce the expressions of the representing
distributions $\eta_x$ and $Z_x$ of the inverse operators
$V^{-1}_k$ and ${}^tV^{-1}_k$ by using the representing
measures~$\mu_x$ and $\nu_x$ of~$V_k$ and~${}^tV_k$. They are
given  by the following formulas {\samepage
\begin{gather*}
\forall\; x \in \mathbb{R}^d,\qquad \eta_x =
{}^tQ(\nu_x),\\ %\label{eq1.13}\\
\forall\; x \in \mathbb{R}^d,\qquad
Z_x = {}^tP(\mu_x),%\label{eq1.14}
\end{gather*}
where ${}^tP$ and ${}^tQ$ are the transposed  operators of $P$ and
$Q$ respectively. }

 The contents of the
paper are as follows. In Section~\ref{sec2} we recall some basic
facts from Dunkl's theory, and describe the Dunkl operators and
the Dunkl kernel. We def\/ine in Section~\ref{sec3} the Dunkl
transform introduced in \cite{5} by C.F.~Dunkl, and we give the
main theorems proved for this transform,  which will be used in
this paper. We study in Section~\ref{sec4} the Dunkl convolution
product and the Dunkl transform of distributions which will be
useful in the sequel, and when the multiplicity function takes
integer values, we give another proof of the geometrical form of
Paley--Wiener--Schwartz theorem for the Dunkl transform. We prove
in Section~\ref{sec5} some inversion formulas for the Dunkl
intertwining operator~$V_k$ and its dual ${}^tV_k$ on spaces of
functions and distributions. Section~\ref{sec6} is devoted to
proving under the
 condition that the
multiplicity function takes integer values an inversion
formula for the Dunkl intertwining operator $V_k$, and we deduce
the expression of the representing distributions of the inverse
operators $V^{-1}_k$ and ${}^tV^{-1}_k$. In Section~\ref{sec7} we
give some applications
 of the preceding
inversion formulas.

\section{The eigenfunction of the Dunkl operators}\label{sec2}

In this section we collect some notations and results on the Dunkl
operators and the Dunkl kernel (see \cite{3, 4, 5, 7, 9, 10}).

\subsection[Reflection groups, root systems and multiplicity
functions]{Ref\/lection groups, root systems and multiplicity
functions}\label{sec2.1}

We consider $\mathbb{R}^d$ with the Euclidean scalar product
$\langle \cdot,\cdot\rangle$ and $\|x\| = \sqrt{\langle x ,
x\rangle}$. On $\mathbb{C}^d$, $\|\cdot\|$ denotes also the
standard Hermitian norm, while $\langle z, w\rangle =
\sum^d_{j=1}z_j \overline{w_j}$ .

For $\alpha \in \mathbb{R}^d \backslash \{0\}$, let
$\sigma_\alpha$ be the ref\/lection in the hyperplane $H_\alpha
\subset \mathbb{R}^d$ orthogonal to $\alpha$, i.e.
\begin{gather*}
\sigma_\alpha(x) = x - \left(\frac{2\langle
\alpha,x\rangle}{\|\alpha\|^2}\right)\alpha .%\label{eq2.1}
\end{gather*}
 A
f\/inite set $R \subset \mathbb{R}^d \backslash \{0\}$ is called a
root  system if $R \cap \mathbb{R}\alpha = \{\pm \alpha\}$ and
$\sigma_\alpha R = R$ for all $\alpha \in R$. For a given  root
system $R$ the ref\/lections $\sigma_\alpha$, $\alpha \in R$,
 generate
a f\/inite group $W \subset O(d)$, the ref\/lection group
associated with $R$. All ref\/lections in $W$ correspond to
suitable pairs of roots. For a given $\beta \in \mathbb{R}^d
\backslash \cup_{\alpha \in R }H_\alpha$, we f\/ix the positive
subsystem $R_+ = \{\alpha \in R ; \langle \alpha, \beta\rangle >
0\}$, then for each $\alpha \in R$ either $\alpha \in R_+$ or
$-\alpha \in R_+$.

A function $k : R \rightarrow \mathbb{C}$ on a root system $R$ is
called a multiplicity function if it is invariant under the action
of the associated ref\/lection group $W$. If one regards $k$ as a
function on the corresponding ref\/lections, this means that $k$
is constant on the conjugacy classes of ref\/lections in $W$. For
abbreviation, we introduce the index
\begin{gather*}
\gamma = \gamma(R) =
 \sum_{\alpha \in R_+}
k(\alpha).%\label{eq2.2}
\end{gather*}
 Moreover, let $\omega_k$ denotes the weight
function
\begin{gather*}
\omega_k(x) = \prod_{\alpha \in R_+}|\langle
\alpha,x\rangle|^{2k(\alpha)},%\label{eq2.3}
\end{gather*} which is $W$-invariant
and homogeneous of degree~$2\gamma$.

For $d = 1$ and $W= \mathbb{Z}_2$, the multiplicity function $k$
is a single parameter denoted $\gamma$ and
\begin{gather*}
\forall\; x \in
\mathbb{R},\qquad \omega_k(x) = |x|^{2\gamma}.%\label{eq2.4}
\end{gather*} We
introduce the Mehta-type constant
\begin{gather*}
c_k = \left( \int_{\mathbb{R}^d}
e^{-\|x\|^2}\omega_k(x)dx\right)^{-1},%\label{eq2.5}
\end{gather*} which is known
for all Coxeter groups $W$ (see~\cite{3, 6}).

\subsection{The Dunkl operators and the Dunkl kernel}\label{sec2.2}

The Dunkl operators $T_j$, $j = 1,\dots,d$, on $\mathbb{R}^d$,
associated with the f\/inite ref\/lection group $W$ and the
multiplicity function~$k$, are given for a function $f$ of class
$C^1$ on $\mathbb{R}^d$ by
\begin{gather*}
T_jf(x) = \frac{\partial}{\partial x_j}f(x) + \sum_{\alpha \in
R_+} k(\alpha) \alpha_j \frac{f(x) -
f(\sigma_\alpha(x))}{\langle\alpha,x\rangle}. %\label{eq2.6}
\end{gather*} In the
case $k \equiv 0$, the $T_j$, $j = 1, 2,\dots,d$, reduce to the
corresponding partial derivatives. In this paper, we will assume
throughout that $k \geq 0$ and $\gamma > 0$.

For $f$ of class $C^1$ on $\mathbb{R}^d$ with compact
\mbox{support} and $g$ of class $C^1$ on $\mathbb{R}^d$ we have
\begin{gather}
\int_{\mathbb{R}^d}T_jf(x) g(x) \omega_k(x)dx =
-\int_{\mathbb{R}^d} f(x) T_j g(x) \omega_k(x) dx,\qquad j = 1,
2,\dots, d.\label{eq2.7}
\end{gather}
 For $y \in \mathbb{R}^d$, the system
\begin{gather}
 T_j u(x, y)  = y_j u(x, y),\qquad j =
1, 2,\dots,d,\nonumber\\ u(0,y)  = 1, \label{eq2.8}
\end{gather}
admits a unique analytic solution on $\mathbb{R}^d$, denoted by
$K(x,y)$ and called the Dunkl kernel.

This kernel has a unique holomorphic extension to $\mathbb{C}^d
\times \mathbb{C}^d$.

\begin{example}\label{example2.1}
From \cite{5}, if $d = 1$ and $W = \mathbb{Z}_2$, the Dunkl kernel
is given by
\begin{gather*}
K(z,t) = j_{\gamma-1/2} (izt) + \frac{zt}{2\gamma+1}
j_{\gamma+1/2} (izt),\qquad z, t \in
\mathbb{C},%\label{eq2.9}
\end{gather*}
 where for $\alpha \geq - 1/2$, $j_\alpha$ is
the normalized Bessel function def\/ined by
\begin{gather*}
j_\alpha(u) = 2^\alpha \Gamma(\alpha + 1)
\frac{J_\alpha(u)}{u^\alpha} = \Gamma(\alpha +1) \sum^\infty_{n=0}
\frac{(-1)^n(u/2)^{2n}}{n!\Gamma(n+\alpha +1)},
\qquad u \in \mathbb{C},%\label{eq2.10}
\end{gather*} with $J_\alpha$  being the Bessel
function of f\/irst kind and index $\alpha$ (see~\cite{16}).
\end{example}

The Dunkl kernel possesses the following properties.
\begin{itemize}\itemsep=0pt
\item[(i)] For $z, t \in \mathbb{C}^d$, we have $K(z,t) = K(t,z)$,
$K(z,0) = 1$, and $K(\lambda z,t) = K(z, \lambda t)$ for all
$\lambda \in \mathbb{C}$.

\item[(ii)] For all $\nu \in \mathbb{Z}^d_+$, $x \in
\mathbb{R}^d$, and $z \in \mathbb{C}^d$ we have
\begin{gather}
|D^\nu_z K(x,z)| \leq \|x\|^{|\nu|} \exp \left[\max_{w \in W}
\langle wx, {\rm Re}\, z\rangle\right].\label{eq2.11}
\end{gather}
 with
\[
D^\nu_z = \frac{\partial^{|\nu|}} {\partial z_1^{\nu_1} \cdots
\partial z^{\nu_d}_d } \qquad \mbox{and} \qquad |\nu| =
 \nu_1 +
\cdots + \nu_d.
\]
 In particular
\begin{gather}
|D^\nu_z K(x,z)|
\leq \|x\|^{|\nu|} \exp [\|x\| \|{\rm Re}\, z\|]],\label{eq2.12} \\
|K(x,z)| \leq \exp [\|x\| \|{\rm Re}\, z\|],\nonumber%\label{eq2.13}
\end{gather}
 and for all
$x, y \in \mathbb{R}^d$
\begin{gather}
|K(ix,y)| \leq 1, \label{eq2.14}
\end{gather}

 \item[(iii)] For all $x, y \in \mathbb{R}^d$
 and $w \in W$ we have
\begin{gather*}
K(-ix, y) = \overline{K(ix, y)} \qquad \mbox{and} \qquad K(wx, wy)
=
  K(x,y).%\label{eq2.15}
  \end{gather*}

 \item[(iv)] The function $K(x,z)$ admits for all $x \in
 \mathbb{R}^d$ and $z \in \mathbb{C}^d$ the following Laplace type
 integral representation
\begin{gather}
K(x,z) = \int_{\mathbb{R}^d} e^{\langle y,z\rangle}d\mu_x(y),
 \label{eq2.16}
\end{gather}
 where $\mu_x$ is the measure given by the relation \eqref{eq1.3} (see~\cite{14}).
\end{itemize}

\begin{remark}\label{remark2.1}
When $d = 1$ and $W = \mathbb{Z}_2$, the relation \eqref{eq2.16}
is of the form
\begin{gather*}
K(x,z) = \frac{\Gamma(\gamma+1/2)}
{\sqrt{\pi}\Gamma(\gamma)}|x|^{-2\gamma}\int^{|x|}_{-|x|}(|x| -
y)^{\gamma-1}(|x| + y)^\gamma e^{yz}dy.%\label{eq2.17}
\end{gather*}
 Then in
this case the measure $\mu_x$ is given for all $x \in
\mathbb{R}\backslash \{0\}$ by $d\mu_x(y) = \mathcal{K}(x,y) dy$
with
\begin{gather*}
\mathcal{K}(x,y) =
\frac{\Gamma(\gamma+1/2)}{\sqrt{\pi}\Gamma(\gamma)}|x|^{-2\gamma}(|x|
-y)^{\gamma-1}(|x|+y)^\gamma 1_{]{-}|x|, |x|[}(y),%\label{eq2.18}
\end{gather*}
where $1_{]{-}|x|, |x|[}$ is the characteristic function of the
interval $]{-}|x|, |x|[$.
\end{remark}

\section{The Dunkl transform}\label{sec3}

In this section we def\/ine the Dunkl transform and we give the
main results satisf\/ied by this transform which will be used in
the following sections (see~\cite{5, 9, 10}).

\smallskip

\noindent
 {\bf Notation.}
We denote by $\mathbb{H}(\mathbb{C}^d)$ the space of entire
functions on $\mathbb{C}^d$ which are rapidly decreasing and of
exponential type. We equip this space with the classical topology.

\smallskip

The Dunkl transform of a function $f$ in
$\mathcal{S}(\mathbb{R}^d)$ is given by
\begin{gather}
\forall\; y \in \mathbb{R}^d,\qquad  \mathcal{F}_D(f)(y) =
\int_{\mathbb{R}^d}f(x)K(x, -iy)\omega_k(x)dx.\label{eq3.1}
\end{gather}
 This
transform satisf\/ies the relation
\begin{gather}
\mathcal{F}_D(f) = \mathcal{F}\circ {}^tV_k(f),\qquad f \in
\mathcal{S}(\mathbb{R}^d),\label{eq3.2}
\end{gather} where $\mathcal{F}$ is the
classical Fourier transform on $\mathbb{R}^d$ given by
\begin{gather*}
\forall\; y \in \mathbb{R}^d,\qquad \mathcal{F}(f)(y) =
\int_{\mathbb{R}^d}f(x) e^{-i\langle x,y\rangle}dx,\qquad f\in
\mathcal{S}(\mathbb{R}^d).%\label{eq3.3}
\end{gather*}
 The following theorems are
proved in \cite{9, 10}.

\begin{theorem} \label{theorem3.1} The transform
$\mathcal{F}_D$ is a topological isomorphism
\begin{itemize}\itemsep=0pt
\item[\rm i)] from $\mathcal{D}(\mathbb{R}^d)$ onto
$\mathbb{H}(\mathbb{C}^d)$, \item[\rm ii)] from
$\mathcal{S}(\mathbb{R}^d)$ onto itself.
\end{itemize}
The inverse transform is given by
\begin{gather}
\forall\; x \in \mathbb{R}^d,\qquad \mathcal{F}^{-1}_D(h)(x) =
\frac{c^2_k}{2^{2\gamma+d}}\int_{\mathbb{R}^d}h(y)K(x,iy)
\omega_k(y)dy.\label{eq3.4}
\end{gather}
\end{theorem}

\begin{remark}\label{remark3.1}
Another  proof of Theorem \ref{theorem3.1} is given in \cite{17}.
\end{remark}

When the multiplicity function satisf\/ies $k(\alpha)\in
\mathbb{N}$ for all $\alpha \in R_+$, M.F.E.~de~Jeu has proved
in~\cite{10} the following geometrical form of Paley--Wiener
theorem for functions.

\begin{theorem}\label{theorem3.2} Let $E$ be a $W$-invariant compact convex set of $
\mathbb{R}^d$ and $f$ an entire function on~$ \mathbb{C}^d$. Then
$f$ is the Dunkl transform of a function in $ \mathcal{D}(
\mathbb{R}^d)$ with support in $E$, if and only if for all $q \in
\mathbb{N}$ there exists a positive constant $C_q$ such that
\begin{gather*}
\forall \; z \in \mathbb{C}^d, \qquad |f(z)| \leq C_q
(1+||z||)^{-q}
e^{I_E({\rm Im}\,z)},%\label{eq3.5}
\end{gather*} where $I_E$ is the
 gauge associated to
the polar of $E$,  given by
\begin{gather}
\forall \; y \in \mathbb{R}^d, \qquad I_E(y) = \sup_{x\in E}
\langle x,y\rangle .\label{eq3.6}
\end{gather}
\end{theorem}

\section{The Dunkl convolution product and the Dunkl transform\\ of distributions}\label{sec4}

\subsection{The Dunkl translation operators and the Dunkl convolution product\\ of functions}\label{sec4.1}

The def\/initions and properties of the Dunkl translation
operators and the Dunkl convolution product of functions presented
in this subsection are given in the seventh section of
\cite[pa\-ges~33--37]{17}.

The Dunkl translation operators $\tau_x$, $x \in \mathbb{R}^d$,
are def\/ined on $\mathcal{E}(\mathbb{R}^d)$
 by
 \begin{gather}
 \forall\; y \in
\mathbb{R}^d, \qquad \tau_x f(y) = (V_k)_x(V_k)_y[V^{-1}_k
(f)(x+y)].\label{eq4.1}
\end{gather}
 For $f$ in $\mathcal{S}(\mathbb{R}^d)$ the
function $\tau_x f$ can also be written in the form
\begin{gather}
\forall\; y \in \mathbb{R}^d, \qquad \tau_x f(y) = (V_k)_x
({}^tV^{-1}_k)_y [{}^tV_k(f)(x+y)].\label{eq4.2}
\end{gather}

Using properties of the operators $V_k$ and ${}^tV_k$ we deduce
that for $f$ in $\mathcal{D}(\mathbb{R}^d)$ (resp.
$\mathcal{S}(\mathbb{R}^d))$ and $x \in \mathbb{R}^d$, the
function $y \rightarrow \tau_x f(y)$ belongs to
$\mathcal{D}(\mathbb{R}^d)$ (resp.~$\mathcal{S}(\mathbb{R}^d))$
and we have
\begin{gather}
\forall\; t \in \mathbb{R}^d,\qquad \mathcal{F}_D(\tau_x f)(t) =
K(ix, t) \mathcal{F}_D(f)(t).\label{eq4.3}
\end{gather}

The Dunkl convolution product of $f$ and $g$ in
$\mathcal{D}(\mathbb{R}^d)$ is the function $f\ast_D g$ def\/ined
by
\begin{gather*}
\forall\; x \in \mathbb{R}^d,\qquad f \ast_D g(x) =
\int_{\mathbb{R}^d}\tau_x f(-y) g(y) \omega_k(y)dy.
%\label{eq4.4}
\end{gather*}
For $f$, $g$ in $\mathcal{D}(\mathbb{R}^d)$
(resp.~$\mathcal{S}(\mathbb{R}^d))$ the function $f\ast_Dg$
belongs to $\mathcal{D}(\mathbb{R}^d)$
(resp.~$\mathcal{S}(\mathbb{R}^d))$ and we have
\begin{gather*}
\forall\; t \in \mathbb{R}^d,\qquad \mathcal{F}_D(f\ast_D g)(t) =
\mathcal{F}_D(f)(t)
\mathcal{F}_D(g)(t).%\label{eq4.5}
\end{gather*}

\subsection{The Dunkl convolution product of tempered  distributions}\label{sec4.2}

\begin{definition}\label{definition4.1}
Let $S$ be in $\mathcal{S}'(\mathbb{R}^d)$ and $\varphi$ in
$\mathcal{S}(\mathbb{R}^d)$. The Dunkl convolution product of $S$
and~$\varphi$  is the function $S\ast_D \varphi$ def\/ined by
\begin{gather*}
\forall\; x \in \mathbb{R}^d,\qquad
 S\ast_D\varphi(x) = \langle S_y,
\tau_x \varphi(-y)\rangle.%\label{eq4.6}
\end{gather*}
\end{definition}

\begin{proposition}\label{proposition4.1}
For $S$ in $\mathcal{S}'(\mathbb{R}^d)$ and $\varphi$ in
$\mathcal{S}(\mathbb{R}^d)$ the function $S \ast_D \varphi$
belongs to $\mathcal{E}(\mathbb{R}^d)$ and we have
\begin{gather*}
T^\mu
(S\ast_D \varphi) = S \ast_D(T^\mu (\varphi)),%\label{eq4.7}
\end{gather*}
 where
\begin{gather*}
T^\mu = T_1^{\mu_1}\circ T_2^{\mu_2} \circ \cdots \circ
T_d^{\mu_d} \qquad \mbox{with} \quad \mu = (\mu_1,
\mu_2,\dots,\mu_d) \in \mathbb{N}^d.
\end{gather*}
\end{proposition}

\begin{proof}  We remark f\/irst that the topology of
$\mathcal{S}(\mathbb{R}^d)$ is also generated by the seminorms
\begin{gather*}
Q_{k,l}(\psi) = \sup_{\substack{|\mu| \leq k\\ x \in
\mathbb{R}^d}}\big(1+||x||^2\big)^{l}|T^{\mu}\psi(x)|, \qquad k,l
\in \mathbb{N}.
\end{gather*}

i) Let $x_0 \in \mathbb{R}^d$. We prove f\/irst that
$S\ast_D\varphi$ is continuous at $x_0$. We have
\[
\forall\; x \in \mathbb{R}^d, \qquad S \ast_D\varphi(x) -
S\ast_D\varphi(x_0) = \langle S_y, (\tau_x \varphi -
\tau_{x_0}\varphi)(-y)\rangle.
\]
 We must
prove that $(\tau_x\varphi - \tau_{x_0}\varphi)$ converges to zero
in $\mathcal{S}(\mathbb{R}^d)$ when $x$ tends to $x_0$.

Let $k, \ell \in \mathbb{N}$ and $\mu \in \mathbb{N}^d$ such that
$|\mu| \leq k$. From \eqref{eq4.3}, Theorem~\ref{theorem3.1} and
the rela\-tions~\eqref{eq2.7},~\eqref{eq2.8} we have
\begin{gather*}
\big(1 +\|y\|^2\big)^\ell T^\mu (\tau_x \varphi - \tau_{x_0}\varphi)(-y)
\\\qquad{} = \frac{i^{|\mu|}c^2_k}{2^{2\gamma+d}}
\int_{\mathbb{R}^d}(1+\|\lambda\|^2)^p K(i\lambda,
-y)(I-\Delta_k)^\ell \Big[\lambda^\mu(K(-ix,\lambda)\\
\qquad{}- K(-ix_0,\lambda)) \mathcal{F}_D(\varphi)(\lambda)\Big]
\frac{\omega_k(\lambda)}{(1+\|\lambda\|^2)^p}d\lambda,
\end{gather*}
with $\lambda^\mu = \lambda_1^{\mu_1}
\lambda_2^{\mu_2}\cdots\lambda^{\mu_d}_d,$ $\Delta_k =
\sum^d_{j=1}T_j^2$  the Dunkl Laplacian and $p \in \mathbb{N}$
such that $p > \gamma + \frac{d}{2}+1$.

Using \eqref{eq2.12} and \eqref{eq2.14} we deduce that
\begin{gather*}
Q_{k,\ell}(\tau_x \varphi - \tau_{x_0}\varphi)
=\sup_{\substack{|\mu| \leq k\\ y \in \mathbb{R}^d}} (1
+\|y\|)^\ell |T^\mu (\tau_x\varphi - \tau_{x_0}\varphi)(-y)|
\rightarrow 0 \qquad \mbox{as} \quad x \rightarrow x_0.
\end{gather*} Then the function $S
*_D \varphi$ is continuous at $x_0$, and thus it is continuous on
$ \mathbb{R}^d$.

 Now we will prove that $S
\ast_D\varphi$ admits a partial derivative on $\mathbb{R}^d$ with
respect to the variable~$x_j$. Let $h \in \mathbb{R}\backslash
\{0\}$. We consider the function $f_h$ def\/ined on $\mathbb{R}^d$
by
\begin{gather*}
f_h(y)
=\frac{1}{h}\big(\tau_{(x_1,\dots,x_j+h,\dots,x_d)}\varphi(-y) -
\tau_{(x_1,\dots,x_j,\dots,x_d)} \varphi(-y)\big) -
\frac{\partial}{\partial x_j} \tau_x \varphi(-y).
\end{gather*}
 Using  the
formula
\begin{gather*}
\forall\; y \in \mathbb{R}^d, f_h(y) = \frac{1}{h}
\int^{x_j+h}_{x_j}\left(\int^{u_j}_{x_j}
\frac{\partial^2}{\partial t^2_j}\tau_{(x_1,\dots,t_j,\dots,x_d)}
\varphi(-y)dt_j\right)du_j,
\end{gather*}
we obtain for all $k, \ell \in \mathbb{N}$ and $\mu \in
\mathbb{N}^d$ such that $|\mu| \leq k$:
\begin{gather}
\forall\; y \in \mathbb{R}^d,\qquad
 (1+\|y\|^2)^\ell
T^\mu f_h(y) \nonumber\\
\qquad\qquad\qquad{}=
\frac{1}{h}\int_{x_j}^{x_j+h}\left(\int^{u_j}_{x_j}\big(1+\|y\|^2\big)^\ell
 T^{\mu}\frac{\partial^2}{\partial
t^2_j}\tau_{(x_1,\dots,t_j,x_d)} \varphi(-y)dt_j\right)du_j.
\label{eq4.8}
\end{gather}
 By applying the preceding method to the function
\begin{gather*}
\big(1+\|y\|^2\big)^\ell T^\mu \frac{\partial^2}{\partial
t^2_j}\tau_{(x_1,\dots,t_j,\dots,x_d)} \varphi(-y),
\end{gather*}
 we deduce
from the relation \eqref{eq4.8} that
\begin{gather*}
Q_{k,\ell}(f_h) = \sup_{\substack{|\mu| \leq k \\ y \in
\mathbb{R}^d}} \big(1+\|y\|^2\big)^\ell |T^\mu f_h(y)| \rightarrow
0\qquad \mbox{as} \quad h \rightarrow 0.
\end{gather*}
 Thus the function $S \ast_D
\varphi(x)$ admits a partial derivative at $x_0$ with respect to
$x_j$ and we have
\begin{gather*}
\frac{\partial}{\partial x_j} S\ast_D \varphi(x_0) = \langle S_y,
\frac{\partial}{\partial x_j}\tau_{x_0}
\varphi(-y)\rangle .%\label{eq4.9}
\end{gather*}
 These results is true on $
\mathbb{R}^d$. Moreover the partial derivatives are continuous on
$ \mathbb{R}^d$. By proceeding in a similar way for partial
derivatives of all order with respect to all variables, we deduce
that $S*_D \varphi$ belongs to~$ \mathcal{E}( \mathbb{R}^d)$.

ii) From the i) we have
\begin{gather*}
\forall \; x \in \mathbb{R}^d, \qquad \frac{\partial}{\partial
x_j} S*_D \varphi (x) = \langle S_y,\frac{\partial}{\partial x_j}
\tau_x \varphi (-y) \rangle.
\end{gather*}
On the other hand using the def\/inition of the Dunkl operator
$T_j$ and the relation
\begin{gather*}
T_j(\tau_x \varphi (-y)) = \tau_x (T_j\varphi) (-y),
\end{gather*}
 we obtain
 \begin{gather*}
  \forall \; x \in \mathbb{R}^d,
\qquad T_j( S*_D \varphi )(x) = \langle S_y, \tau_x (T_j\varphi)
(-y) \rangle  =   S*_D (T_j \varphi)(x).
\end{gather*}
By iteration we get
\begin{gather*}
\forall \; x \in \mathbb{R}^d, \qquad T^\mu( S*_D \varphi )(x) =
S*_D (T^\mu \varphi)(x).\tag*{\qed}
\end{gather*}
\renewcommand{\qed}{}
\end{proof}

\subsection{The Dunkl transform of distributions}\label{sec4.3}

\begin{definition}\label{definition4.2} \qquad {}
 \begin{itemize}\itemsep=0pt
\item[i)] The Dunkl transform of a distribution $S$ in
$\mathcal{S}'(\mathbb{R}^d)$ is def\/ined by
\begin{gather*}
\langle \mathcal{F}_D(S), \psi\rangle = \langle S,
\mathcal{F}_D(\psi)\rangle, \psi \in
\mathcal{S}(\mathbb{R}^d).%\label{eq4.10}
\end{gather*}
\item[ii)] We def\/ine the Dunkl transform of a distribution $S$
in $\mathcal{E}'(\mathbb{R}^d)$ by
\begin{gather}
\forall\; y \in \mathbb{R}^d,\qquad \mathcal{F}_D(S)(y) = \langle
S_x, K(-iy,x)\rangle.\label{eq4.11}
\end{gather}
\end{itemize}
\end{definition}

\begin{remark} \label{remark4.1}
When the distribution $S$ in $\mathcal{E}'(\mathbb{R}^d)$ is given
by the function $g \omega_k$ with $g$ in
$\mathcal{D}(\mathbb{R}^d)$, and denoted by $T_{g \omega_k}$, the
relation \eqref{eq4.11} coincides with \eqref{eq3.1}.
\end{remark}

\noindent {\bf Notation.} We denote by $\mathcal{H}(\mathbb{C}^d)$
the space of entire functions on $\mathbb{C}^d$ which are slowly
increasing and of exponential type.  We equip this space with the
classical topology.

\smallskip

The following theorem is given in \cite[page 27]{17}.

\begin{theorem}\label{theorem4.1}
The transform $\mathcal{F}_D$ is a topological isomorphism
{\samepage
\begin{itemize}\itemsep=0pt
\item[\rm i)] from $\mathcal{S}'(\mathbb{R}^d)$ onto itself; \item[\rm
ii)] from $\mathcal{E}'(\mathbb{R}^d)$ onto
$\mathcal{H}(\mathbb{C}^d)$.
\end{itemize}}
\end{theorem}

\begin{theorem}\label{theorem4.2}
 Let $S$ be in $\mathcal{S}'(\mathbb{R}^d)$ and
$\varphi$ in $\mathcal{S}(\mathbb{R}^d)$. Then,  the distribution
on $\mathbb{R}^d$ given by $(S\ast_D \varphi)\omega_k$ belongs to
$\mathcal{S}'(\mathbb{R}^d)$ and we have
\begin{gather}
\mathcal{F}_D(T_{(S\ast_D \varphi)\omega_k})=
\mathcal{F}_D(\varphi)\mathcal{F}_D(S).\label{eq4.12}
\end{gather}
\end{theorem}

\begin{proof}
i) As $S$ belongs to $\mathcal{S}'(\mathbb{R}^d)$ then there
exists a positive constant $C_0$ and $k_0, \ell_0 \in \mathbb{N}$
such that
\begin{gather}
 |S\ast_D\varphi(x)| = |\langle S_y,
\tau_x\varphi(-y)\rangle|\leq C_0
Q_{k_0,\ell_0}(\tau_x\varphi).\label{eq4.13}
\end{gather}

But by using the inequality
\begin{gather*}
\forall\; x, y \in \mathbb{R}^d,\qquad  1 +\|x+y\|^2 \leq
2\big(1+\|x\|^2\big)\big(1+\|y\|^2\big),
\end{gather*} the relations \eqref{eq4.2},
\eqref{eq1.3} and the properties of the operator
 $^{t}V_k$ (see Theorem~3.2 of~\cite{17}),  we deduce that there exists a positive constant $C_1$ and
$k, \ell \in \mathbb{N}$ such that
\begin{gather*}
Q_{\ell_0,\ell_0}(\tau_x\varphi) \leq
C_1\big(1+\|x\|^2\big)^{\ell_0}Q_{k,\ell}(\varphi).%\label{eq4.14}
\end{gather*}
 Thus
from \eqref{eq4.13} we obtain
\begin{gather}
|S\ast_D \varphi(x)| \leq C\big(1+\|x\|^2\big)^{\ell_0}
Q_{k,\ell}(\varphi), \label{eq4.15}
\end{gather} where
$C$ is a positive constant.  This inequality shows that the
distribution on $\mathbb{R}^d$ associated with the function
$(S*_{D} \varphi)\omega_k$ belongs to
$\mathcal{S}'(\mathbb{R}^d)$.

ii) Let $\psi$ be in $\mathcal{S}(\mathbb{R}^d)$. We shall prove
f\/irst that
\begin{gather}
\langle T_{(S\ast_D \varphi)\omega_k},\psi \rangle = \langle
\check{S}, \varphi \ast_D \check{\psi}\rangle ,\label{eq4.16}
\end{gather}
where $ \check{S}$ is the distribution in $ \mathcal{S'}(
\mathbb{R}^d)$ given by
\begin{gather*}
 \langle \check{S},\phi\rangle = \langle
S, \check{\phi} \rangle,
\end{gather*}
 with
\begin{gather*}
\forall\;x \in \mathbb{R}^d,\qquad \check{ \phi}(x) = \phi(-x).
\end{gather*}
 We consider the two
sequences $\{\varphi_n\}_{n \in \mathbb{N}}$ and $\{\psi_m\}_{m
\in \mathbb{N}}$ in $\mathcal{D}(\mathbb{R}^d)$ which converge
respectively to~$\varphi$ and~$\psi$ in
$\mathcal{S}(\mathbb{R}^d)$. We have
\begin{gather*}
\langle T_{(S\ast_D \varphi_n)\omega_k}, \psi_m\rangle  =  \int_{
\mathbb{R}^d} \langle S_y, \tau_{x}\varphi_n(-y)\rangle \psi_m
(x)\omega_k(x)dx,\\
\phantom{\langle T_{(S\ast_D \varphi_n)\omega_k}, \psi_m\rangle
}{} =  \langle S_y, \int_{\mathbb{R}^d}\!\!\psi_m(x)
\tau_x\varphi_n(-y)\omega_k(x) dx\rangle = \langle S_y,
\int_{\mathbb{R}^d}\!\!\check{\psi}_m(x) \tau_{-x} \varphi_n(-y)
\omega_k(x) dx\rangle.
\end{gather*}

Thus
\begin{gather}
\langle T_{(S\ast_D \varphi_n)\omega_k},\psi_m\rangle =
 \langle \check{S},\varphi_n *_D \check{\psi}_m\rangle.\label{eq4.17}
\end{gather}
 But
\begin{gather*}
\langle T_{(S\ast_D \varphi_n)\omega_k}, \psi_m\rangle -
  \langle T_{(S\ast_D \varphi)\omega_k}, \psi_m\rangle
 = \int_{\mathbb{R}^d} \check{S}\ast_D(\varphi_n - \varphi)(x) \check{\psi}_m(x)
 \omega_k(x)dx.
\end{gather*}
 Thus from \eqref{eq4.15} there exist a positive constant
 $M$ and $k, \ell \in
 \mathbb{N}$ such that
\begin{gather*}
|T_{(S\ast_D \varphi_n)\omega_k},\psi_m\rangle - \langle
T_{(S\ast_D \varphi)\omega_k},
 \psi_m\rangle | \leq M Q_{k,\ell}
 (\varphi_n-\varphi).
\end{gather*}
 Thus
\begin{gather}
\langle T_{(S\ast_D \varphi_n)\omega_k}, \psi_m\rangle
 \underset{n\rightarrow + \infty}{\longrightarrow}
 \langle T_{(S\ast_D \varphi)\omega_k},\psi_m\rangle.\label{eq4.18}
\end{gather}
 On the other hand we have
\begin{gather}
\langle T_{(S\ast_D \varphi)\omega_k},\psi_m\rangle
 \underset{m\rightarrow + \infty}{\longrightarrow}
 \langle T_{(S\ast_D \varphi)\omega_k},\psi
 \rangle,\label{eq4.19}
\end{gather}
 and
\begin{gather}
\varphi_n *_D \check{\psi}_m \underset{\substack{n\rightarrow +
\infty\\ m \rightarrow + \infty}}{\longrightarrow}
 \varphi
 *_D\check{\psi},\label{eq4.20}
 \end{gather}
the limit is in $ \mathcal{S}( \mathbb{R}^d)$.

  We deduce \eqref{eq4.16} from \eqref{eq4.17}, \eqref{eq4.18}, \eqref{eq4.19} and \eqref{eq4.20}.

  We prove now the relation \eqref{eq4.12}. Using \eqref{eq4.16} we obtain for all
  $\psi$ in $\mathcal{S}(\mathbb{R}^d)$
  \begin{gather*}
  \langle \mathcal{F}_D(T_{(S\ast_D \varphi)\omega_k}),\psi\rangle =
  \langle T_{(S\ast_D \varphi)\omega_k}, \mathcal{F}_D(\psi)\rangle,
  = \langle \check{S}, \varphi*_D
  {(\mathcal{F}_D(\psi))}\check{}\rangle.
  \end{gather*}
But
\begin{gather*}
\varphi*_D {(\mathcal{F}_D(\psi))}\check{} =
(\mathcal{F}_D[\mathcal{F}_D(\varphi)\psi])\check{}.
\end{gather*}
 Thus
\begin{gather*}
\langle \breve{S}, \varphi*_D (\mathcal{F}_D(\psi))\check{}\rangle
= \langle S, \mathcal{F}_D[\mathcal{F}_D(\varphi)\psi]\rangle, =
\langle \mathcal{F}_D(\varphi) \mathcal{F}_D(S), \psi\rangle.
\end{gather*}
Then
\begin{gather*}
\langle \mathcal{F}_D(T_{(S\ast_D \varphi)\omega_k}),\psi\rangle =
\langle \mathcal{F}_D(\varphi)\mathcal{F}_D(S),\psi\rangle .
\end{gather*}
 This
completes the proof of \eqref{eq4.12}.
\end{proof}

We consider the positive function $\varphi$ in
$\mathcal{D}(\mathbb{R}^d)$ which is radial for $d \geq 2$ and
even for $d = 1$, with \mbox{support} in the closed ball of center
$0$ and radius 1, satisfying
\begin{gather*}
\int_{\mathbb{R}^d}\varphi(x) \omega_k(x)dx = 1,
\end{gather*}
 and $\phi$
the function on $[0, + \infty[$ given by
\begin{gather*}
\varphi(x) = \phi(\|x\|) = \phi(r) \qquad  \mbox{with}\quad r =
\|x\|.
\end{gather*}
 For $\varepsilon \in ]0,1]$,
we denote by $\varphi_\varepsilon$ the function on $\mathbb{R}^d$
def\/ined by
\begin{gather}
\forall\; x \in \mathbb{R}^d,\qquad \varphi_\varepsilon(x) =
\frac{1}{\varepsilon^{2\gamma+d}}
\phi(\frac{\|x\|}{\varepsilon}).\label{eq4.21}
\end{gather}
This function satisf\/ies the following properties:
\begin{itemize}\itemsep=0pt
\item[i)] Its \mbox{support} is contained in the closed ball
$B_\varepsilon$ of center $0$, and radius $\varepsilon$.
\item[ii)] From \cite[pages 585--586]{13} we have
\begin{gather}
\forall\; y \in \mathbb{R}^d,\qquad
\mathcal{F}_D(\varphi_\varepsilon)(y) = \frac{2^{\gamma +
\frac{d}{2}}}{c_k} \mathcal{F}_B^{\gamma +
\frac{d}{2}-1}(\phi)(\varepsilon\|y\|),\label{eq4.22}
\end{gather} where
$\mathcal{F}_B^{\gamma + \frac{d}{2}-1}(f)(\lambda)$ is the
Fourier--Bessel transform given by
\begin{gather}
\forall\; \lambda \in \mathbb{R}, \qquad \mathcal{F}^{\gamma +
\frac{d}{2}-1}_B(f)(\lambda) = \int^\infty_0f(r)j_{\gamma +
\frac{d}{2}-1}(\lambda r) \frac{r^{2\gamma +d-1}}{2^{\gamma +
\frac{d}{2}}\Gamma\left(\gamma +
\frac{d}{2}\right)}dr,\label{eq4.23}
\end{gather} with $j_{\gamma +
\frac{d}{2}-1}(\lambda r)$ the normalized Bessel function.
\item[iii)] There exists a positive constant $M$ such that
\begin{gather}
\forall\; y \in \mathbb{R}^d,\qquad
|\mathcal{F}_D(\varphi_\varepsilon)(y)-1| \leq \varepsilon
M\|y\|^2.\label{eq4.24}
\end{gather}
\end{itemize}

\begin{theorem}\label{theorem4.3} Let $S$ be in $\mathcal{S}'(\mathbb{R}^d)$. We
have
\begin{gather}
\lim_{\varepsilon \rightarrow 0}(S\ast_D\phi_\varepsilon)\omega_k
= S,\label{eq4.25}
\end{gather} where the
limit is in $\mathcal{S}'(\mathbb{R}^d)$.
 \end{theorem}

\begin{proof}
We deduce \eqref{eq4.25} from \eqref{eq4.12}, \eqref{eq4.22},
\eqref{eq4.24} and Theorem~\ref{theorem4.1}.
\end{proof}

\begin{definition}\label{definition4.3}
Let $S_1$ be in $\mathcal{S}'(\mathbb{R}^d)$ and $S_2$ in
$\mathcal{E}'(\mathbb{R}^d)$. The Dunkl convolution product
of~$S_1$ and~$S_2$ is the distribution $S_1 \ast_D S_2$ on
$\mathbb{R}^d$ def\/ined by
\begin{gather}
\langle S_1 \ast_D S_2,\psi\rangle = \langle S_{1,x}, \langle
S_{2,y}, \tau_x\psi(y) \rangle\rangle,\qquad \psi \in
\mathcal{D}(\mathbb{R}^d).\label{eq4.26}
\end{gather}
\end{definition}

\begin{remark}\label{remark4.2}
The relation \eqref{eq4.26} can also be written in the form
\begin{gather}
\langle S_1 \ast_D S_2,\psi\rangle = \langle S_1, \check{S}_2
\ast_D \psi\rangle.\label{eq4.27}
\end{gather}
\end{remark}

\begin{theorem}\label{theorem4.4} Let $S_1$ be in
$\mathcal{S}'(\mathbb{R}^d)$ and $S_2$ in
$\mathcal{E}'(\mathbb{R}^d)$. Then the distribution $S_1 \ast_D
S_2$ belongs to $\mathcal{S}'(\mathbb{R}^d)$ and we have
\begin{gather*}
\mathcal{F}_D(S_1 \ast_D S_2) =
\mathcal{F}_D(S_2)\cdot \mathcal{F}_D(S_1).%\label{eq4.28}
\end{gather*}
\end{theorem}

\begin{proof}
We deduce the result from \eqref{eq4.27}, the relation
\begin{gather*}
T_{(\check{S_2}*_D \mathcal{F}_D(\psi))\omega_k } = \check{S_2}*_D
T_{\mathcal{F}_D(\psi)\omega_k },
\end{gather*}
 and Theorem~\ref{theorem4.2}.
 \end{proof}

 \subsection[Another proof of the geometrical
  form of the Paley-Wiener-Schwartz theorem
  for the Dunkl transform]{Another proof of the geometrical
  form\\ of the Paley--Wiener--Schwartz theorem
  for the Dunkl transform}\label{sec4.4}

 In this subsection we suppose that the multiplicity
 function satisf\/ies
 $k(\alpha) \in \mathbb{N}
  \backslash \{0\}$ for all $\alpha \in R_+$.

 The main result is to give  another proof of the geometrical form of
Paley--Wiener--Schwartz theorem for the transform $\mathcal{F}_D$,
given in \cite[pages 23--33]{17}.

\begin{theorem}\label{theorem4.5} Let $E$ be a
$W$-invariant compact convex set
 of $\mathbb{R}^d$ and $f$ an entire function on~$\mathbb{C}^d$. Then
 $f$ is the Dunkl transform of a distribution in
 $\mathcal{E}'(\mathbb{R}^d)$ with \mbox{support} in $E$ if and only if
 there exist a positive constant $C$ and $N \in \mathbb{N}$ such
 that
\begin{gather}
\forall\; z \in \mathbb{C}^d, \qquad |f(z)| \leq C (1 + \|z\|^2)^N
  e^{I_E({\rm Im}\, z)},\label{eq4.29}
\end{gather}
  where $I_E$ is the function given by \eqref{eq3.6}.
\end{theorem}

\begin{proof}
\textit{Necessity condition.}
 We consider a distribution $S$ in $\mathcal{E}'(\mathbb{R}^d)$
 with \mbox{support} in $E$.

 Let $\mathcal{X}$ be in $\mathcal{D}(\mathbb{R}^d)$ equal to 1 in
 a neighborhood of $E$, and $\theta$ in $\mathcal{E}(\mathbb{R})$
  such that
\[
\theta(t) = \left\{ \begin{array}{ll}
  1, &\mbox{if} \ t \leq 1,\\
  0, &\mbox{if} \ t > 2.
  \end{array}\right.
  \]
  We put $\eta = {\rm Im}\, z$, $z \in \mathbb{C}^d$ and we take
  $\varepsilon > 0$.
  We denote by $\psi_z$ the function def\/ined on
  $\mathbb{R}^d$ by
\[\psi_z(x) = \chi(x) K(-ix,z)
 |W|^{-1}\sum_{w \in W}
  \theta(\|z\|^\varepsilon (\langle w x,\eta\rangle -
  I_E(\eta))).
  \]
  This function belongs to $\mathcal{D}(\mathbb{R}^d)$
   and as $E$ is $W$-invariant, then it is
  equal to $K(-ix,z)$ in a neighborhood of $E$. Thus
\begin{gather*}
\forall\; z \in \mathbb{C}^d,\qquad \mathcal{F}_D(S)(z)
   = \langle S_x, \psi_z(x)\rangle.%\label{eq4.30}
   \end{gather*}
  As $S$ is with compact \mbox{support}, then it is of f\/inite order $N$.
  Then there exists a positive cons\-tant~$C_0$ such that
\begin{gather}
\forall\; z \in \mathbb{C}^d, \qquad |\mathcal{F}_D(S)(z)| \leq
C_0 \sum_{|p| \leq N}
  \sup_{x \in \mathbb{R}^d}|D^p \psi_z(x)|.\label{eq4.31}
  \end{gather}
  Using the Leibniz rule, we obtain
\begin{gather}
\forall\; x \in \mathbb{R}^d, \qquad D^p \psi_z(x) =
\sum_{q+r+s=p}
  \frac{p!}{q!r!s!} D^q \mathcal{X}
  (x)D^r K(-ix,z)\nonumber\\
\phantom{\forall\; x \in \mathbb{R}^d, \qquad D^p \psi_z(x) =}{}
\times D^s |W|^{-1}
 \sum_{w \in W}\theta(\|z\|^\varepsilon (\langle wx, \eta\rangle
  - I_E(\eta))).\label{eq4.32}
  \end{gather}
  We have
\begin{gather}
\forall\; x \in \mathbb{R}^d, \qquad |D^q\chi(x)| \leq {\rm
const},\label{eq4.33}
\end{gather}
  and if  $M$ is the estimate of $\sup\limits_{t \in \mathbb{R}}
  |\theta^{(k)}(t)|$, $k \leq N$, we obtain
\begin{gather}
\forall\; x \in \mathbb{R}^d, \qquad \left|D^s\left(\sum_{w \in W}
  \theta(\|z\|^\varepsilon (\langle wx,\eta\rangle
  - I_E(\eta)))\right)\right| \leq M(\|z\|^\varepsilon \|\eta\|)^{|s|}.
  \label{eq4.34}
\end{gather}
  On the other hand from \eqref{eq2.11} we have
\begin{gather}
\forall\; x \in \mathbb{R}^d,\qquad
  |D^r K(-ix,z)| \leq \|z\|^r
  e^{\max_{w \in W}\langle wx,\eta\rangle}.
  \label{eq4.35}
\end{gather}
  Using  inequalities \eqref{eq4.33}, \eqref{eq4.34}, \eqref{eq4.35} and \eqref{eq4.32} we deduce
  that there exists a positive cons\-tant~$C_1$ such that
\[
\forall\; x \in \mathbb{R}^d, \qquad |D^p \psi_z(x)|
  \leq C_1
  (1+\|z\|^2)^{N(1+\varepsilon)}
  e^{\max_{w \in W}\langle wx,\eta\rangle}.
\]
  From this relation and \eqref{eq4.31} we obtain
\begin{gather}
\forall\; z \in \mathbb{C}^d,\qquad |\mathcal{F}_D(S)(z)| \leq
  C_2(1+\|z\|^2)^{N(1+\varepsilon)}
  \sup_{x \in E} e^{\max_{w \in W}
  \langle w x,\eta\rangle},\label{eq4.36}
\end{gather}
  where $C_2$ is a positive constant, and
  the supremum is calculated
  when $\|z\| \geq 1$, for
\[
\langle wx, \eta\rangle \leq I_E(\eta) +
  \frac{2}{\|z\|^\varepsilon},
  \]
  because if not we have $\theta = 0$.
  This inequality implies
\begin{gather}
\sup_{x \in E} e^{\max_{w \in W}\langle w x,\eta\rangle}
  \leq e^2 \cdot e^{I_E(\eta)}.\label{eq4.37}
\end{gather}  From \eqref{eq4.36}, \eqref{eq4.37} we deduce that there exists a positive
  constant $C_3$ independent from $\varepsilon$
  such that
\begin{gather*}
\forall\; z \in \mathbb{C}^d, \qquad \|z\| \geq 1, \qquad |
   \mathcal{F}_D(S)(z)|
   \leq C_3(1+\|z\|^2)^{N(1+\varepsilon)}
  e^{I_E(\eta)}.%\label{eq4.38}
\end{gather*}
  If we make $\varepsilon \rightarrow 0$ in this relation we obtain
\eqref{eq4.29} for  $\|z\| \geq 1$. But this inequality is also
true (with
  another constant) for $\|z\| \leq 1$, because in the set $\{z \in \mathbb{C}^d,
  \|z\| \leq 1\}$ the function
  $\mathcal{F}_D(S)(z)e^{- I_E(\eta)}$ is bounded.

\textit{\it Sufficient condition.}
  Let $f$ be an entire function on $\mathbb{C}^d$ satisfying the
  condition~\eqref{eq4.29}. It is clear that the distribution given by
   the restriction of $f\omega_k$ to $\mathbb{R}^d$
  belongs to $\mathcal{S}'(\mathbb{R}^d)$. Thus from Theorem \ref{theorem4.1}i there
  exists a  distribution $S$ in $\mathcal{S}'(\mathbb{R}^d)$ such
  that
\begin{gather}
T_{f\omega_k} = \mathcal{F}_D(S).\label{eq4.39}
\end{gather}
  We shall show that the \mbox{support} of $S$ is contained in $E$. Let
  $\varphi_\varepsilon$ be the function given by the relation \eqref{eq4.21}.
  We consider the distribution
 \begin{gather}
 T_{f_\varepsilon \omega_k}  = \mathcal{F}_D
  (T_{(S \ast_D \varphi_\varepsilon)
   \omega_k}).\label{eq4.40}
  \end{gather}
  From Theorem \ref{theorem4.2} and \eqref{eq4.39}, \eqref{eq4.40} we deduce that
\begin{gather*}
f_\varepsilon = \mathcal{F}_D
  (\varphi_\varepsilon)f.%\label{eq4.41}
\end{gather*}
  The properties of the function $f$ and \eqref{eq4.22}, \eqref{eq4.23} and \eqref{eq4.24}
  show that the function $f_\varepsilon$ can be extended  to an
  entire function on $\mathbb{C}^d$ which satisf\/ies: for all $q \in
  \mathbb{N}$ there exists a positive constant $C_q$ such that
\begin{gather}
\forall\; z \in \mathbb{C}^d,\qquad  |f_\varepsilon(z)| \leq
C_q(1+\|z\|)^{-q}
  e^{I_{E+B_\varepsilon}({\rm Im}\,z)} .\label{eq4.42}
\end{gather}
  Then from \eqref{eq4.42}, Theorem~\ref{theorem3.2} and \eqref{eq4.40},
  the function $(S \ast
  \varphi_\varepsilon)\omega_k$ belongs to
  $\mathcal{D}(\mathbb{R}^d)$ with
  \mbox{support} in $E + B_\varepsilon$.
  But from Theorem~\ref{theorem4.3}, the
  family $(S \ast \varphi_\varepsilon)\omega_k$
  converges
  to $S$ in $\mathcal{S}'(\mathbb{R}^d)$ when $\varepsilon$ tends
  to zero. Thus for all $\varepsilon > 0$, the \mbox{support} of $S$ is
  in $E + B_\varepsilon$,
  then it is contained in~$E$.
  \end{proof}

\begin{remark}\label{remark4.3}
  In the following we give an ameliorated version of the proof of
  Proposition~6.3 of~\cite[page~30]{17}.
  \end{remark}

  Let $E$ be a $W$-invariant compact convex set of $\mathbb{R}^d$
  and $x \in E$. The function $f(x,\cdot)$ def\/ined on $\mathbb{C}^d$ by
\[
f(x,z) = e^{-i\big(\sum\limits^d_{j=1}x_jz_j\big)},
\]
  is entire on $\mathbb{C}^d$ and satisf\/ies
\[
\forall\; z \in \mathbb{C}^d,\qquad |f(x,z)|
  \leq e^{I_E({\rm Im}\, z)}.
  \]
  Thus from Theorem \ref{theorem4.5} there exists a distribution
  $\tilde{\eta}_x$ in $\mathcal{E}'(\mathbb{R}^d)$ with \mbox{support} in
  $E$ such that
\[
\forall\; y \in \mathbb{R}^d, \qquad f(x,y) =
  e^{-i\langle x,y\rangle} = \langle
  \tilde{\eta}_x, K(-iy,\cdot )\rangle.
  \]
  Applying now the remainder of the proof
  given in \cite[page~32]{17}, we
  deduce that the \mbox{support} of the representing distribution
  $\eta_x$ of the inverse Dunkl intertwining operator $V^{-1}_k$
 is contained in~$E$.

 \section{Inversion formulas for the
 Dunkl intertwining operator\\ and its dual}\label{sec5}

 \subsection[The pseudo-differential operators $P$]{The pseudo-dif\/ferential operators $\boldsymbol{P}$}\label{sec5.1}

\begin{definition}\label{definition5.1} We def\/ine the pseudo-dif\/ferential
operator $P$
 on $\mathcal{S}(\mathbb{R}^d)$ by
\begin{gather}
\forall\; x \in \mathbb{R}^d,\qquad P(f)(x) =
 \frac{\pi^d c_k^2}{2^{2\gamma}}
 \mathcal{F}^{-1}[\omega_k \mathcal{F}(f)](x).
 \label{eq5.1}
\end{gather}
\end{definition}

\begin{proposition}\label{proposition5.1} The distribution $T_{\omega_k}$
  given  by the function
 $\omega_k$,
 is in~${\cal S'}(\mathbb{R}^d)$ and for all $f$ in~$\mathcal{S}(\mathbb{R}^d)$ we have
\begin{gather*}
\forall\; x \in \mathbb{R}^d,\qquad P(f)(x) =
  \frac{\pi^d c_k^2}{2^{2\gamma}}
 \mathcal{F}(T_{\omega_k}) * \breve{f}(-x).%\label{eq5.2}
 \end{gather*}
 where $*$ is the classical
 convolution production  of a distribution
 and a function
 on $\mathbb{R}^d$.
 \end{proposition}

\begin{proof}
 It is clear that the distribution $T_{\omega_k}$ given by the
 function $\omega_k$ belongs to ${\cal S'}(\mathbb{R}^d)$. On the
 other hand from the relation \eqref{eq5.1} we have
\begin{gather*}
   \forall\; x \in \mathbb{R}^d,\qquad P(f)(x) =
  \frac{\pi^d c_k^2}{2^{2\gamma}}
 \int_{\mathbb{R}^d}\mathcal{F}(f(\xi+x))(y)\omega_k(y)dy.
\end{gather*}
 Thus
\begin{gather}
\forall\; x \in \mathbb{R}^d,\qquad P(f)(x) =
    \frac{\pi^d c_k^2}{2^{2\gamma}}\langle \mathcal{F}(T_{\omega_k})_y, f(x+y)
\rangle. \label{eq5.3}
\end{gather}
 With the def\/inition of the classical convolution product of a
 distribution and a function on $\mathbb{R}^d$, the relation \eqref{eq5.3}
 can also be written in the form
\begin{gather*}
\forall\; x \in \mathbb{R}^d,\qquad P(f)(x) =
  \frac{\pi^d c_k^2}{2^{2\gamma}}
 \mathcal{F}(T_{\omega_k}) *
 \breve{f}(-x).\tag*{\qed}
 \end{gather*}\renewcommand{\qed}{}
 \end{proof}

 \begin{proposition}\label{proposition5.2} For all $f$ in $\mathcal{S}(\mathbb{R}^d)$
 the function $P(f)$ is of class $C^\infty$ on
 $\mathbb{R}^d$ and we have
\begin{gather}
\forall \; x \in \mathbb{R}^d, \qquad \frac{\partial}{\partial
x_j} P(f)(x) = P\left(\frac{\partial}{\partial
 \xi_j}f\right)(x), \qquad j = 1, 2,\dots,d.\label{eq5.4}
\end{gather}
 \end{proposition}

\begin{proof}
 By derivation under
 the integral sign,
  and by using the relation
\begin{gather*}
\forall\; y \in \mathbb{R}^d,\qquad
  iy_j \mathcal{F}(f)(y) = \mathcal{F}\left(\frac{\partial}
 {\partial \xi_j}f\right)(y), %\label{eq5.5}
 \end{gather*}
  we obtain \eqref{eq5.4}.
 \end{proof}

 \subsection[Inversion formulas for the Dunkl intertwining operator and its dual on the
 space $\mathcal{S}( \mathbb{R}^d)$]{Inversion formulas for the Dunkl intertwining operator\\ and its dual on the
 space $\boldsymbol{\mathcal{S}( \mathbb{R}^d)}$}\label{sec5.2}

 \begin{theorem}\label{theorem5.1} For all $f$ in $\mathcal{S}( \mathbb{R}^d)$ we
 have
\begin{gather}
\forall\; x \in \mathbb{R}^d, \qquad  {}^tV^{-1}_k(f)(x) =
 V_k(P(f))(x).\label{eq5.6}
\end{gather}
 \end{theorem}

 \begin{proof}
 From \cite[Theorem 4.1]{15}  for all $f$ in
 $\mathcal{S}( \mathbb{R}^d)$,  the
 function ${}^tV_k^{-1}(f)$ belongs to $\mathcal{S}( \mathbb{R}^d)$.
 Then from Theorem~\ref{theorem3.1} we have
 \begin{gather}
 \forall\; x \in \mathbb{R}^d,\qquad {}^tV_k^{-1}
 (f)(x) = \frac{c^2_k}{2^{2\gamma+d}}
 \int_{\mathbb{R}^d}
 K(iy,x) \mathcal{F}_D({}^tV_k^{-1}(f))(y) \omega_k(y)dy.
 \label{eq5.7}
 \end{gather}
 But from the relations \eqref{eq3.2}, \eqref{eq1.9}, \eqref{eq1.3}, we have
\[
\forall\; y \in \mathbb{R}^d,\qquad
 \mathcal{F}_D({}^tV^{-1}_k(f))(y) =
 \mathcal{F}(f)(y),
 \]
 and
\[
\forall\; y \in \mathbb{R}^d, \qquad K(iy, x) =
 \mathcal{F}( \breve{\mu}_x)(y),
\]
 where $\breve{\mu}_x$ is the probability
 measure given for a continuous function
  $f$ on $\mathbb{R}^d$ by
\[
\int_{\mathbb{R}^d}f(t)d\check{\mu}_x(t)
  = \int_{\mathbb{R}^d}f(-t)d\mu_x(t).
\]
  Thus \eqref{eq5.7} can also be written in the form
\begin{gather*}
\forall\; x \in \mathbb{R}^d,\qquad
  {}^tV_k^{-1}
  (f)(x) = \frac{c^2_k}{2^{2\gamma+d}}\int_{ \mathbb{R}}
  \mathcal{F}(\breve{\mu}_x)(y) \omega_k(y)
  \mathcal{F}(f)(y)dy.
\end{gather*}
Then by using \eqref{eq5.1},  the properties of the Fourier
transform
  $\mathcal{F}$ and Fubini's theorem we obtain
  \begin{gather*}
  \forall\; x \in \mathbb{R}^d, \qquad {}^tV_k^{-1}(f)(x) =
  \frac{c^2_k}{2^{2\gamma+d}}\int_{\mathbb{R}^d}
  \mathcal{F}[\omega_k \mathcal{F}(f)](y) d\breve{\mu}_x(y) = \int_{\mathbb{R}^d} P(f)(y)d\mu_x(y).
  \end{gather*}
  Thus
 \begin{gather*}
 \forall\; x \in \mathbb{R}^d, \qquad {}^tV_k^{-1}
  (f)(x) = V_k(P(f))(x).\tag*{\qed}
  \end{gather*} \renewcommand{\qed}{}
\end{proof}

\begin{theorem}\label{theorem5.2} For all $f$
  in $\mathcal{S}( \mathbb{R}^d)$ we
  have
\begin{gather}
\forall\; x \in \mathbb{R}^d,\qquad V_k^{-1}(f)(x)
   = P{}^tV_k(f)(x).
  \label{eq5.8}
\end{gather}
\end{theorem}

\begin{proof}
  We deduce the relation \eqref{eq5.8} by replacing $f$ by ${}^tV_k(f)$ in
\eqref{eq5.6} and by using the fact that the operator $V_k$ is an
  isomorphism from $\mathcal{E}( \mathbb{R}^d)$
  onto itself.
\end{proof}

\subsection[Inversion formulas for the dual
Dunkl intertwining operator on the space
$\mathcal{E}'(\mathbb{R}^d)$]{Inversion formulas for the dual
Dunkl intertwining operator\\ on the space
$\boldsymbol{\mathcal{E}'(\mathbb{R}^d)}$}\label{sec5.3}

The dual Dunkl intertwining operator ${}^tV_k$ on $\mathcal{E}'(
\mathbb{R}^d)$ is def\/ined by
\begin{gather*}
\langle {}^tV_k(S), f\rangle = \langle S, V_k(f)\rangle,\qquad
f \in \mathcal{E}( \mathbb{R}^d).%\label{eq5.9}
\end{gather*}

The operator ${}^tV_k$ is a topological isomorphism from
$\mathcal{E}'( \mathbb{R}^d)$ onto itself. The inverse operator is
given by
\begin{gather}
\langle {}^tV^{-1}_k(S),f\rangle = \langle S,
V_k^{-1}(f)\rangle,\qquad f \in \mathcal{E}(
\mathbb{R}^d),\label{eq5.10}
\end{gather}
see \cite[pages 26--27]{17}.

\begin{theorem}\label{theorem5.3} For all $S$ in
$\mathcal{E}'( \mathbb{R}^d)$ the operator ${}^tV_k^{-1}$
satisfies also the relation
\begin{gather}
\langle {}^tV_k^{-1}(S),f\rangle = \langle S, P {}^tV_k(f)\rangle,
\qquad f \in \mathcal{S}( \mathbb{R}^d).\label{eq5.11}
\end{gather}
\end{theorem}

\begin{proof}
We deduce \eqref{eq5.11} from \eqref{eq5.8} and \eqref{eq5.10}.
\end{proof}

 \section{Other expressions
of the inversion formulas\\ for the Dunkl intertwining operator
and its dual\\ when the multiplicity function is
integer}\label{sec6}

In this section we suppose that the multiplicity
 function satisf\/ies
$k(\alpha) \in \mathbb{N}\backslash \{0\}$ for all $\alpha \in
R_+$. The following two Propositions give some other properties of
the operator $P$ def\/ined by \eqref{eq5.1}.

\begin{proposition}\label{proposition6.1}
Let $E$ be a  compact convex set of $\mathbb{R}^d$. Then for all
$f$ in $\mathcal{D}(\mathbb{R}^d)$ we have
\begin{gather*}
\mbox{\rm supp}\, f \subset E \Rightarrow \mbox{\rm supp}\, P(f)
\subset E. %\label{eq6.1}
\end{gather*}
\end{proposition}

\begin{proof}
 From the relation \eqref{eq5.1} we have
 \begin{gather}
  \forall\; x \in \mathbb{R}^d,\qquad P(f)(x)
 = \frac{\pi^d c_k^2}{2^{2\gamma}}
 \int_{\mathbb{R}^d}\mathcal{F}f(y)e^{i\langle x,y
  \rangle}\omega_k(y)dy.\label{eq6.2}
\end{gather}
  We consider the function $F$ def\/ined by
\begin{gather*}
\forall \; z \in \mathbb{C}^d, \qquad F(z) =
  \left(\prod_{ \alpha \in R_+} (\langle\alpha,z\rangle)^{2k(\alpha)}\right)
  \mathcal{F}(f)(z).%\label{eq6.3}
\end{gather*}
This function  is  entire on $\mathbb{C}^d$ and by using
Theorem~2.6 of~\cite{1}  we deduce that for all $q \in
\mathbb{N}$, there exists a positive constant $C_q$ such that
\begin{gather}
\forall\; z \in \mathbb{C}^d, \qquad |F(z)|
 \leq C_q  (1+\|z\|^2)^{-q}e^{I_E
 ({\rm Im}\,z)},
 \label{eq6.4}
\end{gather}
 where $I_E$ is the function given by \eqref{eq3.6}.

 The relation \eqref{eq6.2} can also be written in the form
\begin{gather}
\forall\; x \in \mathbb{R}^d,\qquad P(f)(x) =
 \frac{\pi^d c_k^2}{2^{2\gamma}}
\int_{\mathbb{R}^d}F(y)e^{i\langle x,y\rangle}
 dy.\label{eq6.5}
\end{gather}
Thus \eqref{eq6.5}, \eqref{eq6.4} and Theorem 2.6 of~\cite{1},
imply that $\mbox{\rm supp} \, Pf \subset E$.
\end{proof}

\begin{proposition}\label{proposition} For all $f$ in
 $\mathcal{S}(\mathbb{R}^d)$ we have
\begin{gather}
P(f) = \frac{\pi^d c_k^2}{2^{2\gamma}}\left[ \prod_{\alpha \in R_+
}(-1)^{k(\alpha)}\left(\alpha_1 \frac{\partial}{\partial \xi_1 }
+\cdots+ \alpha_d \frac{\partial}{\partial \xi_d
}\right)^{2k(\alpha)}\right](f). \label{eq6.6}
\end{gather}
\end{proposition}

\begin{proof}  For all $f$ in
$\mathcal{S}(\mathbb{R}^d)$, we have
\begin{gather}
\forall \; y \in \mathbb{R}^d, \quad \omega_k(y) \mathcal{F}(f)(y)
 = \prod_{\alpha \in
R_+}(\langle \alpha,y
\rangle)^{2k(\alpha)}\mathcal{F}(f)(y).\label{eq6.7}
\end{gather}
But
\begin{gather}
\forall\; y \in \mathbb{R}^d, \qquad \langle \alpha,y\rangle
\mathcal{F}(f)(y) = \mathcal{F}\left[- i \left(\alpha_1
\frac{\partial}{\partial \xi_1} + \cdots+ \alpha_d
\frac{\partial}{\partial \xi_d} \right)f\right](y).
 \label{eq6.8}
\end{gather}
 From \eqref{eq6.7}, \eqref{eq6.8}
we obtain
\[
\forall\; y \in \mathbb{R}^d, \qquad \omega_k (y)
\mathcal{F}(f)(y) = \mathcal{F}\left[\prod_{\alpha \in
R_+}(-1)^{k(\alpha)}\left(\alpha_1 \frac{\partial}{\partial \xi_1}
+ \cdots+ \alpha_d \frac{\partial}{\partial
\xi_d}\right)^{2k(\alpha)}f \right](y).
\]

This relation, Def\/inition~\ref{definition5.1} and the inversion
formula for the Fourier transform $\mathcal{F}$ \linebreak
imply~\eqref{eq6.6}.
\end{proof}

\begin{remark}\label{remark6.1}
In this case the operator $P$ is not a pseudo-dif\/ferential
operator but it is a partial dif\/ferential operator.
\end{remark}

\subsection[The
differential-difference operator $Q$]{The
dif\/ferential-dif\/ference operator $\boldsymbol{Q}$}
\label{sec6.1}

\begin{definition}\label{definition6.1}
We def\/ine the dif\/ferential-dif\/ference operator $Q$ on
$\mathcal{S}(\mathbb{R}^d)$ by
\begin{gather*}
\forall\;x \in \mathbb{R}^d,\qquad Q(f)(x) = {}^tV^{-1}_k \circ P
\circ \,{}^t
V_k(f)(x).%\label{eq6.9}
\end{gather*}
\end{definition}

\begin{proposition}\label{proposition6.3} \quad {}
\begin{itemize}\itemsep=0pt
\item[\rm i)] The operator $Q$ is linear and continuous from
$\mathcal{S}(\mathbb{R}^d)$ into itself. \item[\rm ii)] For all
$f$ in $\mathcal{S}(\mathbb{R}^d)$ we have
\[
\forall\;x \in \mathbb{R}^d,\qquad  T_j Q(f)(x) = Q(T_j f)(x),
\qquad j=1,\dots ,d,
\] where
$T_j$, $j = 1, 2,\dots,d$,  are the Dunkl operators.
\end{itemize}
\end{proposition}

\begin{proof}
We deduce the result from the properties of the operator $^{t}V_k$
(see Theorem 3.2 of~\cite{17}),
 and Proposition~\ref{proposition5.2}.
\end{proof}

\begin{proposition}\label{proposition6.4} For all $f$ in
$\mathcal{S}(\mathbb{R}^d)$  we have
\begin{gather}
\forall\;x \in \mathbb{R}^d,\qquad Q(f)(x) = \frac{\pi^d c_k^2
}{2^{2\gamma}}\mathcal{F}_D^{-1}(\omega_k \mathcal{F}_D (f)
)(x).\label{eq6.10}
\end{gather}
\end{proposition}

\begin{proof}
Using the relations \eqref{eq3.2}, \eqref{eq5.1} and the
properties of the operator $^{t}V_k$ (see Theorem~3.2
of~\cite{17}),  we deduce from Def\/inition~\ref{definition6.1}
that  \begin{gather*} \forall\; x \in \mathbb{R}^d,\qquad  Q(f)(x)
= \mathcal{F}_D^{-1}\{\mathcal{F} \circ P ({}^t V_k(f))\}(x) =
\frac{\pi^d c_k^2 }{2^{2\gamma}} \mathcal{F}_D^{-1}\{\mathcal{F}
\circ \mathcal{F}^{-1}[\omega_k \mathcal{F}_D (f)]\}(x).
\end{gather*}
As the function $\omega_k \mathcal{F}_D (f)$  belongs to
$\mathcal{S}(\mathbb{R}^d)$, then by applying the fact that the
classical Fourier transform $\mathcal{F}$ is bijective from $
\mathcal{S}( \mathbb{R}^d)$ onto itself, we obtain
\begin{gather*}
\forall\;x \in \mathbb{R}^d,\qquad Q(f)(x) = \frac{\pi^d c_k^2
}{2^{2\gamma}}\mathcal{F}_D^{-1}(\omega_k
 \mathcal{F}_D (f) )(x).\tag*{\qed}
 \end{gather*}\renewcommand{\qed}{}
\end{proof}

\begin{proposition}\label{proposition6.5} The distribution $T_{\omega_k^2}$
given by the function $\omega_k^2$ is in~$\mathcal{S'}(
\mathbb{R}^d)$ and for all $f $ in~$\mathcal{S} (\mathbb{R}^d)$ we
have
\begin{gather*}
\forall\;x \in \mathbb{R}^d,\qquad  Q(f)(x) = \frac{\pi^d c_k^4
}{2^{4\gamma+d}}\mathcal{F}_D(T_{\omega_k^2})*_D \breve{f}(-x)
,%\label{eq6.11}
\end{gather*} where $*_D$ is the Dunkl convolution product of a
distribution and a function on~$ \mathbb{R}^d$.
\end{proposition}
\begin{proof}
It is clear that the distribution $T_{\omega_k^2}$
 given by the
function ${\omega_k^2}$ belongs to $\mathcal{S'}( \mathbb{R}^d)$.
On the other hand from the relations \eqref{eq6.10}, \eqref{eq3.4}
 and \eqref{eq4.3} we obtain
\begin{gather*}
 \forall\;x \in
\mathbb{R}^d,\qquad Q(f)(x) = \frac{\pi^d c_k^4
}{2^{4\gamma+d}}\int_{
\mathbb{R}^d} \mathcal{F}_D(\tau_x(f))(y) \omega_k^2(y)dy  \nonumber\\
\phantom{\forall\;x \in \mathbb{R}^d,\qquad Q(f)(x)}{}
   =  \frac{\pi^d c_k^4 }{2^{4\gamma+d}} \langle \mathcal{F}
   (T_{\omega_k^2})_y,\tau_x(f)(y) \rangle.%\label{eq6.12}
\end{gather*}
 Thus Def\/inition~\ref{definition4.1} implies
\begin{gather*}
\forall\;x \in \mathbb{R}^d,\qquad Q(f)(x) = \frac{\pi^d c_k^4
}{2^{4\gamma+d}}\mathcal{F}_D(T_{\omega_k^2})*_D
\breve{f}(-x).\tag*{\qed}
\end{gather*}\renewcommand{\qed}{}
\end{proof}

\begin{proposition}\label{proposition6.6} For
all $f$ in $\mathcal{S}( \mathbb{R}^d)$ we have
\begin{gather}
Q(f) = \frac{\pi^d c_k^2}{2^{2\gamma}}\left[\prod_{\alpha \in R_+}
(-1)^{k(\alpha)}(\alpha_1 T_1 + \dots + \alpha_d
T_d)^{2k(\alpha)}\right] (f).\label{eq6.13}
\end{gather}
\end{proposition}

\begin{proof}
For all $f$ in $\mathcal{S}( \mathbb{R}^d)$, we have
\begin{gather}
\forall \; y \in \mathbb{R}^d, \qquad \omega_k(y)\mathcal{F}_D
(f)(y) = \prod_{\alpha \in R_+} (\langle
\alpha,y\rangle)^{2k(\alpha)} \mathcal{F}_D (f)(y).\label{eq6.14}
\end{gather}
 But using \eqref{eq2.7}, \eqref{eq2.8}
  we deduce that
\begin{gather}
\forall \; y \in \mathbb{R}^d, \qquad \langle \alpha,y\rangle
\mathcal{F}_D (f)(y) = \mathcal{F}_D\big[-i (\alpha_1 T_1 + \cdots
+ \alpha_d T_d)f\big](y).\label{eq6.15}
\end{gather}
From \eqref{eq6.14}, \eqref{eq6.15} we obtain
\begin{gather*}
\forall \; y \in \mathbb{R}^d, \qquad \omega_k(y)\mathcal{F}_D
(f)(y) = \mathcal{F}_D \left[\prod_{\alpha \in R_+}
(-1)^{k(\alpha)}(\alpha_1 T_1 + \cdots + \alpha_d
T_d)^{2k(\alpha)}f\right](y).
\end{gather*}
This relation,
Propositions~\ref{proposition6.3},~\ref{proposition6.4} and
Theorem~\ref{theorem3.1} imply \eqref{eq6.13}.
\end{proof}

\subsection{Other expressions of the inversion
 formulas  for the Dunkl intertwining\\ operator and
   its dual on spaces of functions and distributions}\label{sec6.2}

In this subsection we give other expressions of the inversion
formulas for the operators $V_k$ and ${}^tV_k$ and we deduce the
expressions of the representing distributions of the
operators~$V^{-1}_k$ and~${}^tV_k^{-1}$.

\begin{theorem}\label{theorem6.1} For all $f$ in
$\mathcal{S}(\mathbb{R}^d)$ we have
\begin{gather}
\forall\; x \in \mathbb{R}^d,\qquad V^{-1}_k(f)(x) =
{}^tV_k(Q(f))(x). \label{eq6.16}
\end{gather}
\end{theorem}

\begin{proof}
We obtain this result by using of
Proposition~\ref{proposition6.3}, Theorem~\ref{theorem5.2} and
Def\/inition~\ref{definition6.1}.
\end{proof}

\begin{proposition}
\label{proposition6.7} Let $E$ be a $W$-invariant compact convex
set of $\mathbb{R}^d$. Then for all $f$ in
$\mathcal{D}(\mathbb{R}^d)$ we have
\begin{gather}
\mbox{\rm supp}\,   f \subset E \Longleftrightarrow   \mbox{\rm
supp}\,   {}^tV_k(f) \subset E. \label{eq6.17}
\end{gather}
\end{proposition}

\begin{proof}  For all $f$  in
$\mathcal{D}(\mathbb{R}^d)$, we obtain from \eqref{eq3.2} the
relations
\begin{gather*}
^{t}V_k(f) = {\cal F}^{-1}\circ
{\cal F}_D (f),\\ %\label{eq6.18}\\
 ^{t}V_k^{-1}(f) = {\cal F}^{-1}_D\circ
{\cal F} (f).%\label{eq6.19}
\end{gather*}
 We deduce~\eqref{eq6.17} from these relations,
Theorem~\ref{theorem3.2} and Theorem 2.6~of~\cite{1}.
\end{proof}

\begin{proposition}\label{proposition6.8}
Let $E$ be a $W$-invariant compact convex set of $\mathbb{ R}^d$.
Then for all $f$ in $\mathcal{D}(\mathbb{R}^d)$ we have
\begin{gather}
\mbox{\rm supp}\,   f \subset E \Rightarrow   \mbox{\rm supp}\,
Q(f) \subset E. \label{eq6.20}
\end{gather}
\end{proposition}

\begin{proof} We obtain \eqref{eq6.20} from
Def\/inition~\ref{definition6.1},
Propositions~\ref{proposition6.1} and~\ref{proposition6.7}.
\end{proof}

\begin{theorem}\label{theorem6.2} For all $S$ in $\mathcal{E}'( \mathbb{R}^d)$
   the operator ${}^tV^{-1}_k$ satisfies also the relation
\begin{gather}
\langle {}^tV_k^{-1}(S), f\rangle = \langle S,
{}^tV_k(Q(f))\rangle,
   \qquad f \in \mathcal{S}( \mathbb{R}^d).\label{eq6.21}
\end{gather}
\end{theorem}

\begin{proof}   We deduce \eqref{eq6.21} from \eqref{eq5.10} and \eqref{eq6.16}.
\end{proof}

\begin{corollary}\label{corollary6.1} Let $E$ be a $W$-invariant
    compact convex
   set of $\mathbb{R}^d$. For all $S$ in $\mathcal{E}'(
   \mathbb{R}^d)$ with  ${\rm supp}\, S \subset E$, we have
\begin{gather*}
{\rm supp}\,{}^tV_k^{-1}(S) \subset E.%\label{eq6.22}
\end{gather*}
\end{corollary}

\begin{definition}\label{definition6.2} We def\/ine the transposed operators ${}^tP$ and
${}^tQ$ of the operators $P$ and $Q$
on~$\mathcal{S}'(\mathbb{R}^d)$~by
\begin{gather*}
\langle {}^tP(S),f\rangle = \langle S, P(f) \rangle ,\qquad f \in
\mathcal{S}(\mathbb{R}^d),\\ % \label{eq6.23} \\
\langle {}^tQ(S), f\rangle = \langle S,Q(f)\rangle, \qquad f \in
\mathcal{S}(\mathbb{R}^d).%\label{eq6.24}
\end{gather*}
\end{definition}

\begin{proposition}\label{proposition6.9} For all  $S$ in
$\mathcal{S}'(\mathbb{R}^d)$ we have
\begin{gather*}
{}^tP(S) = \frac{\pi^d c_k^2}{2^{2\gamma}} \left[\prod_{\alpha \in
R_+} \left(\alpha \frac{\partial}{\partial \xi_1} + \cdots +
\alpha_d \frac{\partial}{\partial
\xi_d}\right)^{2k(\alpha)}\right]S,\\ %\label{eq6.25}\\
{}^tQ(S) = \frac{\pi^d c_k^2}{2^{2\gamma}} \left[\prod_{\alpha \in
R_+}(\alpha T_1 + \cdots+ \alpha_dT_d)^{2k(\alpha)}
\right]S,%\label{eq6.26}
\end{gather*}
 where $T_j$, $j = 1, 2,\dots,d$, are the
Dunkl operators defined on $\mathcal{S}'(\mathbb{R}^d)$ by
\[
\langle T_j S,f\rangle = - \langle S, T_jf\rangle,\qquad f\in
\mathcal{S}(\mathbb{R}^d).
\]
\end{proposition}

\begin{proposition}\label{proposition6.10} For all $S$ in
$\mathcal{S'}(\mathbb{R}^d) $ we have
\begin{gather*}
\mathcal{F}^{-1}({}^tP(S)) = \frac{\pi^d c_k^2}{2^{2\gamma}}
\mathcal{F}^{-1}(S)
\omega_k,\\ %\label{eq6.27}\\
 \mathcal{F}^{-1}_D({}^tQ(S)) =
\frac{\pi^d c_k^2}{2^{2\gamma}} \mathcal{F}^{-1}_D(S)
\omega_k.%\label{eq6.28}
\end{gather*}
\end{proposition}

\begin{proof}
 We deduce these relations from~\eqref{eq5.1}, \eqref{eq6.10} and
 the def\/initions of the classical Fourier transform
and the Dunkl transform of tempered distributions  on
$\mathbb{R}^d$.
\end{proof}

\begin{theorem}\label{theorem6.3} The representing distributions $\eta_x$ and $Z_x$
of the inverse of the Dunkl intertwining operator and  its dual,
are given by
\begin{gather}
\forall\; x \in \mathbb{R}^d, \qquad \eta_x =
{}^tQ(\nu_x)\label{eq6.29}
\end{gather}
 and
\begin{gather}
\forall\;x \in \mathbb{R}^d, \qquad Z_x =
{}^tP(\mu_x),\label{eq6.30}
\end{gather}
where $\mu_x$ and $\nu_x$ are the representing measures of the
Dunkl intertwining opera\-tor~$V_k$ and
 its dual~${}^tV_k$.
 \end{theorem}

\begin{proof}
From \eqref{eq1.6}, for all $f$ in $\mathcal{S}(\mathbb{R}^d)$ we
have
\begin{gather}
\forall\; x \in \mathbb{R}^d,\qquad  {}^tV_k(Q(f))(x) =
\langle\nu_x, Q(f)\rangle = \langle
{}^tQ(\nu_x),f\rangle.\label{eq6.31}
\end{gather} On the
other hand from \eqref{eq1.4}
\[
\forall\; x \in \mathbb{R}^d,\qquad V^{-1}_k(f)(x) = \langle
\eta_x,f\rangle.
\]
 We obtain \eqref{eq6.29} from
this relation, \eqref{eq6.31} and \eqref{eq6.16}.

Using \eqref{eq1.3}, for all $f$ in $\mathcal{S}(\mathbb{R}^d)$ we
can also write the relation \eqref{eq5.6} in the form
\begin{gather}
\forall\; x \in \mathbb{R}^d, \qquad {}^tV^{-1}_k(f)(x) = \langle
\mu_x, P(f)\rangle =
\langle{}^tP(\mu_x),f\rangle.\label{eq6.32}
\end{gather}
 But from \eqref{eq1.8} we have
\[
\forall\; x \in \mathbb{R}^d, \qquad {}^tV^{-1}(f)(x) = \langle
Z_x,f\rangle.
\]
 We deduce \eqref{eq6.30} from this relation and
\eqref{eq6.32}.
\end{proof}

\begin{corollary}\label{corollary6.2} We have
\begin{gather*}
\forall\;
 x \in \mathbb{R}^d, \qquad \eta_x =
\frac{\pi^d c_k^2}{2^{2\gamma}} \left[\prod_{\alpha \in R_+
}\left(\alpha_1T_1+\cdots+ \alpha_d T_d
\right)^{2k(\alpha)}\right](\nu_x)%\label{eq6.33}
\end{gather*}
 and
\begin{gather*}
\forall\; x \in \mathbb{R}^d, \qquad Z_x = \frac{\pi^d
c_k^2}{2^{2\gamma}}\left[\prod_{\alpha \in R_+}\left(\alpha_1
\frac{\partial}{\partial \xi_1}+\cdots+ \alpha_d
\frac{\partial}{\partial
\xi_d}\right)^{2k(\alpha)}\right](\mu_x).%\label{eq6.34}
\end{gather*}
\end{corollary}

\begin{proof}
We deduce these relations from Theorem~\ref{theorem6.3}
 and Proposition~\ref{proposition6.9}.
 \end{proof}

\section{Applications}\label{sec7}

\subsection[Other proof of the sufficiency
 condition of Theorem~\ref{theorem4.4}]{Other proof of the suf\/f\/iciency
 condition of Theorem~\ref{theorem4.4}}\label{sec7.1}

 Let $f$ be an entire function on $\mathbb{C}^d$ satisfying the
 condition \eqref{eq4.29}. Then from Theorem~2.6 of~\cite{1}, the
 distribution $\mathcal{F}^{-1}(f)$ belongs to $\mathcal{E}'(
 \mathbb{R}^d)$ and we have
\[
{\rm supp}\, \mathcal{F}^{-1}(f) \subset E.
\]

From the relation
\[
\mathcal{F}^{-1}_D(f) = {}^tV_k^{-1} \circ \mathcal{F}^{-1}(f)
\]
 given in \cite[page~27]{17} and Corollary~\ref{corollary6.1}, we deduce that the
 distribution $\mathcal{F}_D^{-1}(f)$ is in $\mathcal{E}'(
 \mathbb{R}^d)$ and its  support is contained in $E$.

 \subsection{Other expressions of the  Dunkl  translation operators}\label{sec7.2}

 We consider the  Dunkl  translation operators
 $\tau_x$, $x \in \mathbb{R}^d$, given by
 the rela\-tions~\eqref{eq4.1},~\eqref{eq4.2}.

 \begin{theorem}\label{theorem 7.1} \quad {}
  \begin{itemize}\itemsep=0pt
 \item[\rm i)] When the multiplicity function
 $k(\alpha)$ satisfies $k(\alpha) >
 0$  for all $\alpha  \in R_+$, we have
\begin{gather}
\forall\; y \in \mathbb{R}^d,\qquad
  \tau_x(f)(y) = \mu_x \ast
 \mu_y(P {}^tV_k(f)),\qquad f \in
  \mathcal{S}( \mathbb{R}^d),\label{eq7.1}
  \end{gather}
 where $\ast$ is the classical convolution product of measures on~$\mathbb{R}^d$.
 \item[\rm ii)] When the multiplicity function satisfies
  $k(\alpha)\in \mathbb{N}\backslash \{0\}$ for all
 $\alpha \in R_+$,  we have
\begin{gather}
\forall\;   y \in \mathbb{R}^d,\qquad \tau_x(f)(y) =
 \mu_x\ast \mu_y({}^tV_k(Q(f))),\qquad f \in \mathcal{S}
 ( \mathbb{R}^d).\label{eq7.2}
 \end{gather}
 \end{itemize}
\end{theorem}

\begin{proof}
i) From the relations \eqref{eq4.1} and \eqref{eq1.3},
 for $f$ in $\mathcal{S}(\mathbb{R}^d)$
  we have
\[
\forall\;x ,y \in \mathbb{R}^d, \qquad \tau_x(f)(y) =
\int_{\mathbb{R}^d}\int_{\mathbb{R}^d}
 V^{-1}_k(f)(\xi+\eta)d\mu_x(\xi)d\mu_y(\eta).
 \]
  By using the def\/inition of the classical convolution product of two
 measures with compact support on $\mathbb{R}^d$, we obtain
\[
\forall\; x,y \in \mathbb{R}^d, \qquad \tau_x(f)(y) = \mu_x \ast
\mu_y(V^{-1}_k(f)).
\]
 Thus Theorem~\ref{theorem5.2} implies the relation~\eqref{eq7.1}.

ii) The same proof as for i)  and Theorem~\ref{theorem6.1} give
the
 relation~\eqref{eq7.2}.
 \end{proof}

\subsection*{Acknowledgements}
  The author would like to thank the referees for their interesting and useful
  remarks.

\pdfbookmark[1]{References}{ref}
\LastPageEnding

 \end{document}